\documentclass[10pt]{amsart}
\usepackage{amsmath}
\usepackage{amsxtra}
\usepackage{amscd}
\usepackage{amsthm}
\usepackage{amsfonts}
\usepackage{amssymb}
\usepackage{eucal}

\theoremstyle{definition}
\theoremstyle{remark}

\newcommand{\nc}{\newcommand}
\nc{\renc}{\renewcommand}
\nc{\ssec}{\subsection}
\nc{\sssec}{\subsubsection}
\nc{\on}{\operatorname}
\nc{\wt}{\widetilde}
\nc{\ol}{\overline}

\nc{\CM}{{\mathcal M}}
\nc{\CN}{{\mathcal N}}
\nc{\CF}{{\mathcal F}}
\nc{\CQ}{{\mathcal Q}}
\nc{\CG}{{\mathcal G}}
\nc{\CY}{{\mathcal Y}}
\nc{\CC}{{\mathcal C}}
\nc{\CO}{{\mathcal O}}
\nc{\CL}{{\mathcal L}}
\nc{\CA}{{\mathcal A}}
\nc{\CP}{{\mathcal P}}
\nc{\CZ}{{\mathcal Z}}
\nc{\CJ}{{\mathcal J}}

\nc{\cg}{{\check\fg}}
\nc{\cG}{{\check G}}
\nc{\cB}{{\check B}}
\nc{\cU}{{\check U}}
\nc{\cH}{{\check H}}

\nc{\BC}{{\mathbb C}}
\nc{\BK}{{\mathbb K}}
\nc{\BZ}{{\mathbb Z}}

\nc{\bC}{{\mathbf C}}
\nc{\bD}{{\mathbf D}}
\nc{\bB}{{\mathbf B}}

\nc{\bg}{{\mathbf g}}
\nc{\bh}{{\mathbf h}}

\nc{\fg}{{\mathfrak g}}
\nc{\fd}{{\mathfrak d}}
\nc{\fD}{{\mathfrak D}}
\nc{\fh}{{\mathfrak h}}
\nc{\fp}{{\mathfrak p}}
\nc{\fq}{{\mathfrak q}}

\nc{\cNt}{\check\CN^{th}}
\nc{\CNt}{\CN^{th}}
\nc{\tg}{\wt{\fg}}
\nc{\tN}{\wt{\CN}}
\nc{\tNt}{\wt{\CN}^{th}}
\nc{\tcN}{\wt{\check\CN}}
\nc{\tcNt}{\wt{\check\CN}^{th}}

\nc{\sF}{{\mathsf{F}}}
\nc{\Hom}{\on{Hom}}
\nc{\Vect}{\on{Vect}}
\nc{\End}{\on{End}}
\nc{\QCoh}{\on{QCoh}}
\nc{\Coh}{\on{Coh}}
\nc{\Bun}{\on{Bun}}
\nc{\Ind}{\on{Ind}}
\nc{\act}{\on{act}}
\nc{\Pro}{\on{Pro}}
\nc{\Rep}{\on{Rep}}
\nc{\St}{\on{St}}
\nc{\one}{{\mathbf{1}}}
\nc{\GUb}{{\overline{G/U}}}
\nc{\cGUb}{{\overline{\check{G}/\check{U}}}}
\nc{\uV}{{\underline{V}}}
\nc{\ul}{{\underline{\ell}}}
\nc{\Fl}{\on{Fl}}
\nc{\tFl}{\wt{\on{Fl}}}
\nc{\Gr}{\on{Gr}}
\nc{\perf}{\on{perf}}
\nc{\bHo}{{\mathbf {Ho}}}
\nc{\wXi}{\wt{\Xi}}

\nc\Spec{\on{Spec}}
\nc\Sch{\on{Sch}}
\renc{\mod}{\text{-mod}}

\title{The notion of category over an algebraic stack}

\author{Dennis Gaitsgory}

\address{Department of Mathematics, Harvard University, 1 Oxford street, Cambridge MA, 02138}

\email{gaitsgde@math.harvard.edu}

\begin{document}

\begin{abstract}
The goal of this note is to spell out the (apparently well-known and intuitively clear)
notion of an abelian category over a stack. In the future we will discuss the (much less
evident) notion, when instead of an abelian category one considers a triangulated one.
\end{abstract}

\maketitle

\date{July, 2005}

\noindent{\bf 1.}
Let $\CC$ be a $\BC$-linear abelian category. We will assume that $\CC$ is closed under 
inductive limits, \footnote{In what follows by an inductive limit we will mean a limit
taken over a small filtering category} i.e., that the tautological embedding $\CC\to \Ind(\CC)$ 
admits a right adjoint $limInd:\Ind(\CC)\to \CC$, and that the latter functor is exact. 
In particular, it makes  sense to tensor objects of $\CC$ by vector spaces.

\medskip

\noindent{\bf 2. The affine case.}
Let $A$ be a commutative algebra. We say that $\CC$ is $A$-linear if we are given
a map $A\to Z(\CC)$, i.e., if $A$ acts functorially on every $\Hom(X,Y)$ for $X,Y\in \CC$.
We shall also say that in this case $\CC$ "lives over $S=\Spec(A)$".

\medskip

We claim that we have a well-defined functor of tensor product
$$M,X\mapsto M\underset{A}\otimes X:A\mod\times \CC\to \CC:$$

If $M=A^I$ for some index set $I$, then $M\underset{A}\otimes X:=X^I$, and
if $M=\on{coker}(A^I\to A^J)$, then 
$$M\underset{A}\otimes X:=\on{coker}(X^I\to X^J),$$
where the $(i,j)$-entry of the matrix $X^I\to X^J$ is given by the action
of the $(i,j)$-entry of the matrix $A^I\to A^J$.

\medskip

\noindent{\bf Lemma 3.}
{\it The above definition is independent of the presentation of $M$ as a quotient.}

\medskip

By construction, the functor of tensor product commutes with inductive limits with respect
to both $M$ and $X$, and is right-exact. In addition, we have:

\medskip

\noindent{\bf Lemma 4.}
{\it 
\smallskip

\noindent {\em (a)} If $M$ is flat, then the functor
$X\mapsto M\underset{A}\otimes X$ is exact.

\smallskip

\noindent {\em (b)}
If $0\to M_1\to M_2\to M_3\to 0$ is a short exact sequence of $A$-modules with $M_3$
flat, then the sequence
$$0\to M_1\underset{A}\otimes X\to M_2\underset{A}\otimes X\to M_3\underset{A}\otimes X\to 0$$
is also short exact.

\smallskip

\noindent {\em (c)} If $M$ is projective and finitely generated, and $M^\vee$ is the dual
module, then the above functor admits left and right adjoints, both given by
$X\mapsto M^\vee\underset{A}\otimes X$.}

\medskip

\begin{proof}

First, let us note that if $M$ is a projective module given by an idempotent of
$A^I$ for some set $I$, then $M\underset{A}\otimes X$ is given by the
corresponding idempotent of $M^I$. This implies that the functor of
tensor product with a projective $A$-module is exact. This implies
point (a), since every flat $A$-module can be represented as an inductive limit 
of projective ones.

Similarly, for point (b) we can assume that $M_3$ is projective,
in which case the short exact sequence splits and the assertion is obvious.

Point (c) is immediate, since we have the adjunctions maps
$$X\to M^\vee \underset{A}\otimes (M\underset{A}\otimes X)\simeq
M\underset{A}\otimes (M^\vee \underset{A}\otimes X) \text{ and }
M^\vee \underset{A}\otimes (M\underset{A}\otimes X)\simeq
M\underset{A}\otimes (M^\vee \underset{A}\otimes X) \to X$$
that satisfy the necessary conditions.

\end{proof}

Finally, we have:

\medskip

\noindent{\bf Proposition 5.}
{\it Assume that $A'$ is a faithfully-flat algebra over $A$. Then 
$A'\underset{A}\otimes X\neq 0$ if $X\neq 0$.}

\medskip

\begin{proof}(Drinfeld)

\medskip

\noindent{\bf Lemma 6.}
{\it If $A'$ is a faithfully flat algebra over $A$, then the quotient $A'/A$ is $A$-flat.}

\medskip

Clearly, the lemma implies Proposition 5, by Lemma 4(b). 

\end{proof}

\begin{proof}(of the Lemma)

It is enough to show that $A'/A\underset{A}\otimes A'$ is $A'$-flat. But
$$A'/A\underset{A}\otimes A'\simeq \on{coker}(A'\to A'\underset{A}\otimes A'),$$
and the latter is a split injection.

\end{proof}

We shall say that $X\in \CC$ is flat over $A$ (or $S$) if the functor 
$M\mapsto M\underset{A}\otimes X:A\mod\to \CC$ is exact. 

\medskip

\noindent{\bf 7. Change of rings.} Let $f:\Spec(A')=S'\to S=\Spec(A)$ is a 
morphism of affine schemes, corresponding to a homomorphism of algebras $A\to A'$.
There exists a universal $A'$-linear category $\CC'$, which admits an
$A$-linear functor $\CC\to \CC'$. We will denote this category by
$\CC\underset{S}\times S'$, and it is constructed as follows:

Objects of $\CC'$ are objects $X\in \CC$, endowed with an additional
action of $A'$, such that the two actions of $A$ (one coming from $A\to A'$, 
and another from $A\to \End(\CM)$), coincide. Morphisms in $\CC'$
are arrows $X_1\to X_2$ in $\CC$ that commute with the the $A'$-action.

The functor $\CC\to \CC'$ is given by $X\mapsto A'\underset{A}\otimes X$,
and it will be denoted by $f^*$. This functor is the left adjoint to the
forgetful functor $f_*:\CC'\to \CC$.

\medskip

Set $S''=S'\underset{S}\times S'$, and let $\CC''$ denote the corresponding base-changed 
category over $S''$. One naturally defines the category of descent data on $\CC'$ with respect 
to $S''$. We will denote it by $Desc_{S''}(\CC')$, and we have a natural functor 
$\CC\to Desc_{S''}(\CC')$.

\medskip

\noindent{\bf Proposition 8.}
{\it Suppose that $S'$ is faithfully-flat over $S$. Then
$\CC\to Desc_{S''}(\CC')$ is an equivalence.}

\medskip

\begin{proof}

This is proved by the usual argument, using Lemma 5. 

\end{proof}

\medskip

\noindent{\bf 9. Stacks: approach I.}    
Let $\CY$ be a stack (algebraic in the faithfully flat sense),
for which the diagonal morphism $\CY\to \CY\times \CY$ is affine.
This is equivalent to demanding that any morphism $S\to \CY$, with $S$ an affine scheme,
is affine. We are going to introduce the notion of sheaf of abelian categories over $\CY$.
In particular, we will obtain a notion of category over a separated scheme.

\medskip

Let $\Sch_\CY^{aff}$ be the category of affine schemes over $\CY$, endowed with
the faithfully flat topology. A sheaf of categories $\CC^{sh}$ over $\CY$ is the following data:

\begin{itemize}

\item
For each $S=\Spec(A)\in \Sch_\CY^{aff}$, a category $\CC_S$ over $S$.

\item
For $f:S_2\to S_1\in \Sch_\CY^{aff}$, an $S_1$-linear functor $f^*:\CC_{S_1}\to \CC_{S_2}$, which
induces an equivalence $\CC_{S_1}\underset{S_1}\times S_2\to \CC_{S_2}$.

\item
For two morphisms $S_3\overset{g}\to S_2\overset{f}\to S_1 \in \Sch_\CY^{aff}$
an isomorphism of functors $g^*\circ f^*\simeq (f\circ g)^*$, such that the natural
compatibility axiom for 3-fold compositions holds.

\end{itemize}

\medskip

Given a sheaf of categories $\CC^{sh}$ over $\CY$ one can form a single category, denoted
$\Gamma(\CY,\CC^{sh})$ or $\CC_\CY$ (or simply $\CC$, where no confusion is likely
to occur) as follows: 

Let $S\to \CY$ be a faithfully flat cover. We define the category $\CC_\CY$ to be 
the category of descent data of $\CC_S$ with respect to the two maps 
$S\underset{\CY}\times S\rightrightarrows S$. Proposition 8 insures that $\CC_\CY$ 
is well-defined, i.e., is canonically independent of the choice of the cover $S$. 

\medskip

Again, by Proposition 8, we have the natural functor
$X\mapsto X_S:\CC_\CY\to \CC_S$ for any $S\in \Sch_\CY^{aff}$, and for
$f:S_2\to S_1$ a functorial isomorphism $f^*(X_{S_1})\simeq X_{S_2}$.

\medskip

When $\CY$ is itself an affine scheme $S=\Spec(A)$, a data of a sheaf of categories $\CC^{sh}$
over $S$ is equivalent to a single category over $S$, which is reconstructed as 
$\CC_S$. In this case we will often abuse the notation and not distinguish
between $\CC^{sh}$ and $\CC_S$.

\medskip

We will now define a functor 
$$\CF,X\mapsto \CF\ast X:\QCoh_\CY\times \CC_\CY\to \CC_\CY.$$

Let $\CF$ be a quasi-coherent sheaf of $\CY$; for $S=\Spec(A)\in \Sch_\CY^{aff}$
we will denote by $\CF_S$ the corresponding quasi-coherent sheaf of $S$. For
$X\in \CC_\CY$ we define
$$(\CF\ast X)_S:=\CF_S\underset{A}\otimes X_S,$$
which by descent gives rise to an object of $\CC_\CY$.

The above functor has the following properties:

\begin{itemize}

\item(i)
$\QCoh_\CY\times \CC_\CY\to \CC_\CY$ is right exact 
and commutes with inductive limits.

\item(ii)
We have a functorial isomorphism
$\CO_\CY\ast X\simeq X$.

\item(iii)
We have functorial isomorphisms
$\CF_1\ast (\CF_2\ast X)\simeq (\CF_1\underset{\CO_\CY}\otimes \CF_2)\ast X$,
compatible with triple tensor products and the isomorphism of (ii).

\end{itemize}

By construction, the assertions of Lemma 4 hold in the present
context, when we replace $M\underset{A}\otimes X$ by $\CF\mapsto \CF\ast X$.

\medskip

\noindent{\bf 10. Descent of categories.} 
Let $f:\CY'\to \CY$ be a map of stacks, and $\CC^{sh}$ a sheaf of categories over $\CY$.
It is clear that it gives rise to a sheaf of categories $\CC'{}^{sh}:=
\CC^{sh}\underset\CY\times \CY'$ over $\CY'$, 
such that for $S\in \Sch_{\CY'}^{aff}$ the category $\CC'_S$ is by definition $\CC_S$,
where $S$ is regarded as an object of $\Sch_{\CY}^{aff}$.

If $g:\CY''\to \CY'$, it is clear that we have an equivalence of sheaves of categories 
$$\CC^{sh}\underset{\CY}\times \CY''\simeq (\CC^{sh}\underset\CY\times \CY')
\underset{\CY'}\times \CY''.$$

\medskip

Suppose now $\CC'{}^{sh}$ is a sheaf of categories over $\CY'$. Let $p^i_j$ be the projection 
on the $j$-th factor from the $i$-fold Cartesian product $\CY^{(i)}$ of $\CY'$ 
over $\CY$ to $\CY'$. Let $\CC^{(i)}_j{}^{sh}$ denote the corresponding base-changed
sheaf of categories categories over $\CY^{(i)}$.

Suppose we are give an equivalence of sheaves of categories over $\CY^{(2)}:
\CC^{(2)}_1{}^{sh}\simeq \CC^{(2)}_2{}^{sh}$; a natural transformation between the two functors
$\CC^{(3)}_1{}^{sh}\to \CC^{(3)}_3{}^{sh}$, such that the two natural transformations between 
the two functors $\CC^{(4)}_1{}^{sh}\to \CC^{(4)}_4{}^{sh}$ coincide.

\medskip

\noindent{\bf Proposition 11.}
{\it Suppose that $\CY'$ is faithfully flat over $\CY$. Then there exists a well-defined
sheaf of categories $\CC^{sh}$ over $\CY$ with an equivalence 
$\CC'{}^{sh}\simeq \CC^{sh}\underset\CY\times \CY'$, and
which gives rise to the above functors and natural transformations.}

\medskip

\begin{proof}

The assertion readily reduces to the case when both $\CY$ and $\CY'$ are affine schemes,
$\Spec(A)$ and $\Spec(A')$, respectively, Let $\Phi$ denote the functor
$\CC^{(2)}_1\to \CC^{(2)}_2$, and $T$ the natural transformation between 
the functors $\Phi^{1,3}$ and $\Phi^{2,3}\circ \Phi^{1,2}$ between $\CC^{(3)}_1$ and
$\CC^{(3)}_3$.

We define $\CC$ to have as objects $X'\in \CC'$ endowed with an isomorphism
$$\alpha_{X'}:\Phi((p_2^1)^*(X'))\to (p_2^2)^*(X'),$$ such that the diagram
$$
\CD
\Phi^{1,3}((p_3^1)^*(X')) @>{T}>>  \Phi^{2,3}\circ \Phi^{1,2}((p_3^1)^*(X')) \\
@VVV   @VVV \\
(p_3^3)^*(X') @<<< \Phi^{2,3}((p_3^2)^*(X'))
\endCD
$$
commutes. Morhisms in this category are $\CC'$-morphisms, commuting with the data
of $\alpha_{X'}$. Evidently, this is an $A$-linear category.

By construction, we have a functor $\CC\to \CC'$, which gives rise to a functor
$$\CC\underset{\Spec(A)}\times \Spec(A')\to \CC'.$$
The fact that the latter is an equivalence is shown by the base-change technique
as in the context of quasi-coherent sheaves.

\end{proof}

\medskip

\noindent{\bf 12. Example: categories with a group-action.}
Let us consider an example of the above situation, when $\CY'=\on{pt}$, $\CY=\on{pt}/G$, 
where $G$ is an affine algebraic group. Let $\CC'{}^{sh}$ be a sheaf of categories over
$\CY'$, i.e., a plain category. Then the data of an equivalence $\CC^{(2)}_1\to \CC^{(2)}_2$ 
together with a natural transformation as above is what can be reasonably called
an action of the group $G$ on $\CC'$.

\medskip

Let us spell this notion out in more detail. We claim that an action of $G$ on a category
category $\CC'$ is equivalent to a data of a functor
$$\on{act}^*:\CC'\to \CO_G\mod\otimes \CC',$$
(here $A\mod\otimes \CC'$ denotes the same thing as $\CC'\underset{\on{pt}}\times \Spec(A)$),
and two functorial isomorphisms related to this functor. This first isomorphism is between the 
identity functor on $\CC'$ and the composition 
$\CC'\overset{\on{act}^*}\to \CO_G\mod\otimes \CC'\to \CC'$,
where the second arrow corresponds to the restriction to $1\in G$.

\medskip

To formulate the second isomorphism, note that from the existing data
we obtain a natural functor $$\on{act}^*_A: A\mod\otimes \CC'\to
\CO_G\mod\otimes A\mod \otimes \CC'\simeq (\CO_G\otimes A)\mod\otimes \CC'$$
for any algebra $A$.

The second isomorphism is between the two functors $\CC\to \CO_{G\times G}\mod\otimes \CC$
that correspond to the two circuits of the diagram
$$
\CD
\CC'  @>{\on{act}^*}>> \CO_G\mod\otimes \CC'  \\
@V{\on{act}^*}VV   @V{\on{act}^*_{\CO_G}}VV   \\
\CO_G\mod\otimes \CC' @>{\on{mult}^*}>> \CO_{G\times G}\mod\otimes \CC',
\endCD
$$
where $\on{mult}$ denoted the multiplication map $G\times G\to G$. These functors
must satisfy the usual compatibility conditions.

\medskip

From Proposition 11, it follows that an action of $G$ on a category $\CC'$ is
equivalent to the data of a sheaf of categories $\CC^{sh}$ over $\on{pt}/G$. 
(As we shall see later, the latter can be also reformulated as a category 
with an action of the tensor category $\Rep(G)$.)

By definition, $\CC:=\CC_{\on{pt}/G}$ can be reconstructed as the category of $G$-equivariant 
objects of $\CC'$. By definition, the latter consists of $X'\in \CC'$, endowed with an isomorphism
$\alpha_X:\on{act}^*(X')\simeq \CO_G\otimes X'$, which is compatible with unit and associativity
constraints. Morphisms in the category are $\CC'$-morphisms, compatible with the
data of $\alpha$.

\medskip

\noindent{\bf 13. Example: categories acted on by a groupoid.}
Generalizing the above set-up, let $S$ be a base scheme, and $\CG\overset{p_2}{\underset{p_1}\rightrightarrows} S$ be an affine groupoid, such
that the maps $p_1,p_2$ (or, equivalently, one of them) are flat. Let $\CC^{sh}$ be a sheaf 
of categories over the quotient stack $\CY=S/\CG$. This data can be rewritten as a 
sheaf of categories $\CC'{}^{sh}$ over $S$, acted on by $\CG$, which means the following:
\footnote{To simplify the notation, we will assume here that $S$ is affine as well.}

\medskip

We must be given a functor
$\act^*:\CC'\to \CC'\underset{S,p_1}\times \CG$,
which is $\CO_S$-linear if we regard $\CC'\underset{S,p_1}\times \CG$ as a category over $S$ via
$\CG\overset{p_2}\to S$, and two functorial isomorphisms related to it. The first isomorphism
is a unit constraint, i.e., an isomorphism between the functor
$$\CC'\overset{\act^*}\to \CC'\underset{S,p_1}\times \CG
\overset{1_\CG^*}\to \left(\CC'\underset{S,p_1}\times \CG\right)\underset{\CG,1_\CG}\times S
\simeq \CC',$$
where $1_\CG:S\to \CG$ is the unit map. 

Formulate the second isomorphism note that for any scheme $S'$, mapping to $S$, 
we obtain a functor
$$\act^*\times \on{id}_{S'}:\CC'\underset{S}\times S'\to
\CC'\underset{S,p_1}\times (\CG\underset{p_2,S}\times S').$$
The second isomorphism is an associativity constraint, i.e., an 
isomorphism between the two functors in the diagram
$$
\CD
\CC'  @>{\act^*}>> \CC'\underset{S,p_1}\times \CG \\
@V{\act^*}VV    @V{\act^*\times \on{id}_\CG}VV   \\
\CC'\underset{S,p_1}\times \CG  @>{\on{mult}^*}>> 
\CC'\underset{S,p_1}\times \left(\CG\underset{p_2,S,p_1}\times \CG\right)
\endCD
$$
such that the natural compatibility conditions hold.

\medskip

\noindent{\bf Lemma 14.}
{\it 
\smallskip

\noindent{\em (a)}
Let $0\to \CF_1\to \CF_2\to \CF_3\to 0$ be a short exact sequence of quasi-coherent
sheaves on $\CG$ with $\CF_3$ being $\CO_S$-flat with respect to $p_2$. Then for
$X\in \CC'$, the sequence
$$0\to \CF_1\ast \act^*(X)\to
\CF_2\ast \act^*(X)\to 
\CF_3\ast \act^*(X)\to 0$$
is also short exact.

\smallskip

\noindent{\em (b)}
If $X\in \CC'$ is $\CO_S$-flat, then $\act^*(X)$ is $\CO_\CG$-flat.}

\medskip

\begin{proof}

Let $S'$ be a scheme with a map $\phi:S'\to \CG$; let $\psi_i=p_i\circ \phi$,
$i=1,2$. We claim that there exists a natural $\CO_{S'}$-linear equivalence
$$\act^*_\phi:\CC'\underset{S,\psi_2}\times S'\to \CC'\underset{S,\psi_1}\times S',$$
defined by
$$X\mapsto (\phi\times \on{id}_{S'})^*\circ (\act^*\times \on{id}_{S'})(X),$$
where $\phi\times \on{id}_{S'}: S'\to \CG\underset{p_2,S,\psi_2}\times S'.$
Its quasi-inverse is defined using the map $\gamma\circ \phi:S'\to \CG$,
where $\gamma$ is the inversion on $\CG$.

We apply this to $S'=\CG$ and $\phi=\gamma$. We obtain an equivalence
$$\act^*_\gamma:\CC'\underset{S,p_1}\times \CG\to \CC'\underset{S,p_2}\times \CG,$$
such that for $X\in \CC'$,
$$\act_\gamma^*(\act^*(X))\simeq p_2^*(X).$$

This readily implies both points of the lemma.

\end{proof}

We say that an object $X\in \CC'$ is $\CG$-equivariant, if we are given
an isomorphism
$$p_1^*(X)\simeq \act^*(X)\in \CC'\underset{S,p_1}\times \CG,$$
compatible with the unit and associativity constraints. Let us denote by
$\CC'{}^{\CG}$ the category of $\CG$-equivariant objects in $\CC'$.

From the definitions we obtain:

\medskip

\noindent{\bf Lemma 15.}
{\it 
\smallskip

\noindent{\em (a)}
For any $X\in \CC'$, the object $(p_1)_*(\act^*(X))$ is naturally $\CG$-equivariant.

\smallskip

\noindent{\em (b)}
The functor $X\mapsto (p_1)_*(\act^*(X))$ is the right adjoint to the 
forgetful functor $\CC'{}^{\CG}\to \CC'$.}

\medskip

In addition, we have:

\medskip

\noindent{\bf Lemma 16.}
{\it Assume that $\CG$ is flat {\em over} $S\times S$. Then every $\CG$-equivariant
object of $\CC'$ is $\CO_S$-flat.}

\medskip

\begin{proof}

This follows from the fact that for $\CF\in \QCoh_S$ and $X\in \CC'$,
$$\act^*(\CF\ast X)\simeq p_2^*(\CF)\ast \act^*(X).$$

\end{proof}

\medskip

\noindent{\bf 17. Stacks: approach II.}   
Let $\Vect_\CY$ denote the tensor category of locally free sheaves of finite
rank on $\CY$. 

Assume now that we are given a category $\CC_\CY$ endowed with an action of the
tensor category $\Vect_\CY$:
$$\ast:\Vect_\CY\times \CC_\CY\to \CC_\CY,$$
which is exact. I.e., for a fixed $\CP\in \Vect_\CY$ the functor $X\mapsto \CP\star X$
is exact, and whenever
$0\to \CP_1\to \CP_2\to \CP_3\to 0$ is
a short exact sequence of objects of $\Vect_\CY$, the corresponding 
sequence
$$0\to \CP_1\ast X\to \CP_2\ast X\to \CP_3\ast X\to 0$$
is also exact. We shall call such a data "a category over $\CY$".

\medskip 

We will now make an additional assumption on the stack $\CY$:

\begin{itemize}

\item The stack $\CY$ is locally Noetherian and every quasi-coherent sheaf on it
is an inductive limit of coherent ones.

\item Every coherent sheaf on $\CY$ can be covered by an object
of $\Vect_\CY$.

\end{itemize}

As in the affine case, this implies that every flat quasi-coherent sheaf on $\CY$
can be represented as an inductive limit of objects of $\Vect_\CY$.

\medskip

\noindent{\bf Theorem 18.}
{\it Under the above assumption on $\CY$, a data of a category over $\CY$ is equivalent
to that of a sheaf of categories over $\CY$.}

\medskip

The rest of this subsection and the next one are devoted to the proof of this theorem. One
direction has been explained above: given a sheaf of categories $\CC^{sh}$ over $\CY$, we
reconstruct $\CC_\CY$ as $\Gamma(\CY,\CC^{sh})$. To carry out the construction in the opposite
direction we will use the above additional assumption on $\CY$. 

\medskip

We claim that the above data extends to an action of the monoidal category $\QCoh_\CY$ on
$\CC_\CY$, satisfying the conditions (i),(ii),(iii) of Sect. 9 and assertions (a), (b) and (c)
of Lemma 4.

First we define an action of the monoidal category $\Coh_\CY$ on $\CC_\CY$:
By assumption, every $\CF\in \Coh_\CY$ can be represented as
$\on{coker}(\CP\to \CQ)$ with $\CP,\CQ\in \Vect_\CY$. We set
$$\CF\ast X: =\on{coker}(\CP\ast X\to \CQ\ast X).$$

To show that this is well-defined, we must consider a commutative diagram
of objects of $\Vect_\CY$
$$
\CD
0 & & 0 & & 0 \\
@AAA   @AAA  @AAA \\
\CP @>{\phi}>> \CQ @>>> \CF @>>> 0\\
@AAA  @AAA  @A{\on{id}}AA & & \\
\CP' @>{\phi'}>> \CQ' @>>> \CF @>>> 0\\
@AAA  @AAA @AAA \\
\CP'' @>{\phi''}>>  \CQ'' @>>> 0
\endCD
$$
with exact rows and columns, and show that the map
$$\on{coker}(\phi\ast \on{id}_X)\to \on{coker}(\phi'\ast \on{id}_X)$$
is an isomorphism.  But this follows from the assumption.

It is clear that the resulting functor is right-exact and satisfies 
properties (ii) and (iii) of Sect. 9. 

\medskip

Next, we have to extend the above action of $\Coh_\CY$ on $\CC_\CY$
to that of $\QCoh_\CY$ by setting for $\CF\simeq \underset{\longrightarrow}{lim}\, \CF_i$
with $\CF_i\in \Coh_\CY$, $\CF\in \QCoh_\CY$,
$$\CF\ast X:=\underset{\longrightarrow}{lim}\, \CF_i\ast X.$$
The resulting action satisfies properties 
(i), (ii) and (iii) of Sect. 9. By assumption, the functor of tensor
product with an object of $\Vect_\CY$ is exact. This implies properties
(a) and (b) Lemma 4, by repeating the proof of {\it loc.cit.}
Property (c) stated in Lemma 4 is evident.

\medskip

\noindent{\bf 19.}
We shall now show how the data of an action of $\QCoh_\CY$ on $\CC_\CY$ with the
above properties reconstructs the categories $\CC_S$ for $S\in \Sch_\CY^{aff}$.

Let $S=\Spec(A)$, let us denote by $\CO_S$ the direct image
of the structure sheaf of $S$ onto $\CY$, regarded as an algebra
in $\QCoh_\CY$. We introduce $\CC_S$ as the category, consisting of objects $X$
of $\CC_\CY$, endowed with an associative action $\CO_S\ast X\to X$, and the morphisms
being $\CC_\CY$-morphisms compatible with the action. 

We have a map of algebras $A\otimes \CO_\CY\to \CO_S$ in $\QCoh_\CY$; this makes 
$\CO_S$ into an $A$-linear category. We also have a functor
$\CC_\CY\to \CC_S$ given by $X\mapsto X_S:=\CO_S\ast X$.

\medskip

Let $f:S_2=\Spec(A_2)\to \Spec(A_1)=S_1$ be a morphism in $\Sch_\CY^{aff}$.
We define a functor $f^*:\CC_{S_1} \to \CC_{S_2}$ by
$$X\mapsto \CO_{S_2}\underset{\CO_{S_1}}\ast X,$$
where for an algebra $\CA$ in $\QCoh_\CY$, a sheaf $\CF$ of $\CA$-modules
and an object $X\in \CC_\CY$ acted on by $\CA$, we set
$$\CF\underset{\CA}\ast X:=\on{coker}((\CA\underset{\CO_\CY} \CA) \ast X 
\rightrightarrows \CF\ast X).$$

We claim that the induced functor $\CC_{S_1}\underset{S_1}\times S_2\to \CC_{S_2}$
is an equivalence. This follows from the fact that $\CO_{S_2}\simeq A_2\underset{A_1}\otimes \CO_{S_1}$. 

Note in addition that for $X\in \CC_\CY$,
\begin{equation*}  \label{one}
\CO_{S_2}\ast X\simeq A_2\underset{A_1}\otimes (\CO_{S_1}\ast X).
\end{equation*}

This implies that for $X\in \CC_\CY$, we have a natural isomorphism 
$f^*(X_{S_1})\simeq X_{S_2}$.

Thus, we have constructed a sheaf of categories over $\CY$, and it remains to show that
the initial category $\CC_\CY$ can be reconstructed by the descent procedure. The usual proof
for coherent sheaves works, once we establish the following:

\medskip

\noindent{\bf  Lemma 20.}
{\it If $\CA\in \QCoh_\CY$ be an algebra, faithfully flat over $\CO_\CY$.
Then the functor $X\mapsto \CA\ast X:\CC_\CY\to \CC_\CY$ is exact and faithful.}

\medskip

\begin{proof}

The exactness part follows by property (a) of Lemma 4.
The faithfulness part follows as in Proposition 5 using 
property (b) of Lemma 4.

\end{proof}

\medskip

\noindent{\bf 21. Example: de-equivariantization.}   
Let $\CY$ be the stack $\on{pt}/G$, where $G$ is an affine algebraic group. Given
a category $\CC$, a structure of category over $\on{pt}/G$ on it is by definition the
same as an action of the tensor category $\Rep(G)$ of 
finite-dimensional representations of $G$ on it:
$$V\in \Rep(G), X\in \CC\mapsto V\ast X\in \CC,$$
which has the exactness property of Sect. 17.

By Theorem 18, such a data gives rise to a sheaf of categories $\CC^{sh}$ over
$\on{pt}/G$.

\medskip

Let us show how to reconstruct the category $\CC':=\CC^{sh}\underset{\on{pt}/G}\times \on{pt}$.
By definition, this is the category, whose objects are
$X'\in \CC$, endowed with an associative action $\CO_G\ast X'\to X'$, and morphisms
are $\CC$-morphisms, compatible with this action.

\medskip

According to \cite{AG}, this data can be rewritten as follows. For every $V\in \Rep(G)$ 
we must be given a map
$$\beta_V:V\ast X\to  X\otimes \underline{V},$$
for every $V\in \Rep(G)$ (here $\underline{V}$ denoted the vector space underlying
a representation), which satisfy the compatibility conditions of \cite{AG}, Sect. 2.2.
One easily shows that the maps $\beta_V$ are necessarily isomorphisms.
Morphisms in this category are $\CC$-morphisms, compatible with the data
of $\beta$.

Thus, $\CC'$ is the category of Hecke eigen-objects in $\CC$ with respect to the action of
$\Rep(G)$. By construction, $\CC'$ carries a canonical action of $G$. Explicitly,
for $X'\in \CC'$, the $\CO_G$-family $\act^*(X')$ is isomorphic to $X'\otimes \CO_G$
as an object of $\CC$. The isomorphisms $\beta$ are given by
$$V\ast (X'\otimes \CO_G)\overset{\beta_V\otimes \on{id}_{\CO_G}}\longrightarrow
X' \otimes \underline{V}\otimes \CO_G \to X' \otimes \underline{V}\otimes \CO_G,$$
where the second arrow is given by the co-action map $\underline{V}\to
\underline{V}\otimes \CO_G$.

According to Sect. 6, the category $\CC$ is reconstructed from $\CC'$ as the
category $\CC'{}^G$ of $G$-equivariant objects. We will refer to $\CC'$ as the de-equivariantization
of $\CC$.

\medskip

\noindent{\bf 22. Another example}   
Generalizing the previous example, let us take $\CY=S/G$, where $S=\Spec(A)$ is an affine
scheme, and $G$ an affine algebraic group acting on it. Let $\CC$ be an abelian category.
A structure on $\CC$ of category over $S/G$ is by definition expressed as follows: 

An action of the tensor category $\Rep(G)$ on $\CC$:
$V,X\mapsto V\ast X$ as above,
and a functorial map $\alpha_X:A\ast X\to X$, where $A$ is regarded as an algebra
in $\Rep(G)$, such that for $V\in \Rep(G)$ the diagram
$$
\CD
V\ast (A\ast X) @>{\sim}>> A\ast (V\ast X) \\
@V{V\ast \alpha_X}VV  @V{\alpha_{V\ast X}}VV  \\
V\ast X @>{\on{id}}>> V\ast X
\endCD
$$
commutes.

\end{document}